\RequirePackage{ifpdf}
\ifpdf 
\documentclass[pdftex]{sigma}
\else
\documentclass{sigma}
\fi

\begin{document}

\allowdisplaybreaks

\renewcommand{\PaperNumber}{042}

\FirstPageHeading

\ShortArticleName{Local Quasitriangular Hopf Algebras}

\ArticleName{Local Quasitriangular Hopf Algebras}

\Author{Shouchuan ZHANG~$^{\dag\,\ddag}$, Mark D. GOULD~$^\ddag$ and Yao-Zhong ZHANG~$^\ddag$}

\AuthorNameForHeading{S.~Zhang, M.D.~Gould and Y.-Z.~Zhang}

\Address{$^\dag$~Department  of Mathematics, Hunan University, Changsha  410082, P.R. China} 
\EmailD{\href{mailto:z9491@yahoo.com.cn}{z9491@yahoo.com.cn}}

\Address{$^\ddag$~Department of Mathematics, University of Queensland,
Brisbane 4072, Australia}
\EmailD{\href{mailto:mdg@maths.uq.edu.au}{mdg@maths.uq.edu.au}, \href{mailto:yzz@maths.uq.edu.au}{yzz@maths.uq.edu.au}}

\ArticleDates{Received January 31, 2008, in f\/inal form April 30,
2008; Published online May 09, 2008}

\Abstract{We f\/ind a new class of Hopf algebras,  local quasitriangular Hopf
algebras, which  generalize quasitriangular Hopf algebras. Using
these Hopf algebras, we obtain solutions of the Yang--Baxter equation
in a systematic way. The category of  modules with f\/inite cycles
over a local quasitriangular Hopf algebra is a braided tensor category.}


\Classification{16W30; 16G10}

\section{Introduction}

The Yang--Baxter equation f\/irst came up in the paper by Yang as
factorization condition of the scattering S-matrix in the many-body
problem in one dimension and in the work by Baxter on exactly
solvable models in statistical mechanics. It has been playing an
important role in mathematics and physics (see \cite{BD82,YG89}). Attempts to f\/ind solutions of the Yang--Baxter equation in
a systematic way have led to the theory of quantum groups and
quasitriangular Hopf algebras (see \cite{Dr86,Ja96}).

 Since the category of  modules with f\/inite
cycles over a local quasitriangular Hopf algebra is a~braided tensor
category, we may also f\/ind solutions of the Yang--Baxter equation in
a~systematic way.

The main results in this paper are summarized in the following
statement.

\begin{theorem}\label{summarized}
(i) Assume that  $(H, \{R_n \})$ is a local quasitriangular
Hopf algebra. Then   $(_H {\cal M}^{\rm cf},$ $C^{\{R_n\}})$, $(_H {\cal
M}^{\rm dcf}, C^{\{R_n\}})$ and  $(_H {\cal M}^{\rm df}, C^{\{R_n\}})$ are
braided tensor
 categories. Furthermore,  if $(M, \alpha ^-)$ is an  $H$-module with finite
cycles and $R_{n+1} = R_n +W_n$ with $W_n \in H_{n+1} \otimes
H_{(n+1)}$,  then $(M, \alpha ^-, \delta ^- )$ is a Yetter--Drinfeld
$H$-module.

(ii) Assume that $B$ is a finite dimensional Hopf algebra and $M$ is
a finite dimensional $B$-Hopf bimodule.  Then $((T_{B} (M)) ^{\rm cop}
\bowtie _\tau T_{B^*} ^c(M^*), \{ R_n\})$ is a local quasitriangular
Hopf algebra.  Furthermore, if $(kQ^a, kQ^c)$ and  $(kQ^s, kQ^{sc})$
 are arrow dual pairings  with finite Hopf quiver $Q$, then
both  $((kQ^a)^{\rm cop} \bowtie _\tau kQ^c, \{ R_n\})$ and
  $((kQ^{s})^{\rm cop} \bowtie _\tau kQ^{sc}, \{ R_n\})$ are  local
quasitriangular Hopf algebras.
\end{theorem}

\section{Preliminaries}\label{s0}

Throughout, we work over a f\/ixed f\/ield $k$. All algebras,
coalgebras, Hopf algebras, and so on, are def\/ined over $k$. Books
\cite{DNR01,Mo93,Sw69,Lu93} provide the necessary background for Hopf
algebras and book \cite {ARS95} provides a nice description of the
path algebra approach.

 Let  $V$ and $W$ be two vector spaces.  $\sigma _V$  denotes the map
 from $V$ to $V^{**}$ by def\/ining $\langle \sigma _V(x), f\rangle = \langle f, x\rangle$
for any $f \in V^*$, $x \in V$. $C_{V,W}$ denotes the map from
$V\otimes  W$ to $W \otimes V$ by def\/ining $C_{V,W}(x\otimes y) =y
\otimes x$ for any $x\in  V$, $y \in W$.  Denote $P$ by $\sum P'
\otimes P''$ for  $P\in V\otimes W$. If $V$ is a f\/inite-dimensional
vector space over f\/ield $k$ with $V^* = {\rm Hom}_k (V, k)$. Def\/ine maps
$b_V: k \rightarrow V\otimes V^*$ and $d_V: V ^* \otimes V
\rightarrow k $ by
\[
b_V(1) = \sum _iv_i \otimes v_i^*  \qquad \mbox{and}\qquad
\sum _{i, j } d_V (v_i^* \otimes v_j ) =\langle v_i^*, v_j \rangle,
\] where $\{
v_i \mid i = 1, 2, \dots , n \}$  is any basis of $V$ and $\{v_i^*
\mid i = 1, 2, \dots , n \}$ is its dual basis in $V^*$. $d_V$  and
$b_V$ are called   evaluation and coevaluation of $V$, respectively.
It is clear  $(d_V \otimes id _U) (id _U \otimes b_V) =id_ U$  and $
(id_V \otimes d _V)(b_V \otimes id _V)= id_ V.$ $\xi _V$ denotes the
linear isomorphism from $V$ to $V^*$ by sending~$v_i$ to~$v_i^*$ for
$i = 1, 2, \dots , n$.   Note that we can def\/ine evaluation $d_V$
when $V$ is inf\/inite.

We will use $\mu$ to denote the multiplication map of an algebra and
use $\Delta$ to denote the comultiplication of a coalgebra. For a
(left or right) module and a (left or right) comodule, denote by
$\alpha ^-$, $\alpha ^+$, $\delta ^-$ and $\delta ^+$ the left
module, right module, left comodule and right comodule structure
maps, respectively. The Sweedler's sigma notations for coalgebras
and comodules are $\Delta (x) = \sum x_1\otimes x_2$, $\delta ^-
(x)= \sum x_{(-1)} \otimes x_{(0)}$, $\delta ^+ (x)= \sum x_{(0)}
\otimes x_{(1)}$. Let $(H, \mu, \eta , \Delta, \epsilon )$ be a~bialgebra and let
 $\Delta ^{\rm cop} := C_{H,H}\Delta$  and $\mu^{\rm op} := \mu C_{H,H}.$ We
denote $(H, \mu, \eta,  \Delta ^{\rm cop}, \epsilon)$  by $H^{\rm cop}$ and
$(H, \mu^{\rm op}, \eta, \Delta, \epsilon )$  by $H^{\rm op}$. Sometimes, we
also denote the unit element of $H$ by $1_H.$

 Let $A$ and $H$ be two bialgebras with
$\varnothing \not= X \subseteq A$ ,   $\varnothing \not= Y \subseteq H$
and $P \in Y \otimes X$, $R\in Y \otimes Y. $ Assume that  $\tau $ is
a linear map from $H\otimes A$ to $k.$   We give the following
notations.

$(Y,R)$  is called almost cocommutative if the following condition
satisf\/ied:
\[
{\rm (ACO)}: \ \ \sum y_2R' \otimes y_1R'' = \sum R'y_1 \otimes R''y_2 \quad
\mbox{for any} \quad y\in Y.
\]

$\tau $ is called a skew pairing on $H\otimes A$ if for any $x, u\in
H$, $y, z \in A$ the following conditions are satisf\/ied:
\begin{gather*}
{\rm (SP1)}: \ \ \tau (x, yz ) = \sum \tau (x_1, y) \tau (x_2, z);\\
{\rm (SP2)}: \ \ \tau (xu, z ) = \sum \tau (x, z_2) \tau (u, z_1);\\
{\rm (SP3)}: \ \ \tau (x ,  \eta ) = \epsilon _H(x); \\
{\rm (SP4)}: \ \ \tau (\eta , y ) = \epsilon _A (y).
\end{gather*}

$P$ is called a copairing of  $Y\otimes X$ if for any $x, u\in H$, $y,
z \in A$ the following conditions are satisf\/ied:
\begin{gather*}
{\rm (CP1)}: \ \ \sum P' \otimes P''_1 \otimes P''_2 = \sum P'Q' \otimes Q''
\otimes P'' \quad \mbox{with}\quad P=Q;\\
{\rm (CP1)}: \ \ \sum P'_1 \otimes P'_2 \otimes P'' = \sum P' \otimes Q'
\otimes P''Q'' \quad \mbox{with} \quad P=Q;\\
{\rm (CP3)}: \ \ \sum P'\otimes \epsilon _A(P'') = \eta _H; \\
{\rm (CP4)}: \ \ \sum \epsilon _H(P') \otimes P'' = \eta _A.
\end{gather*}

For $R \in H \otimes H$ and two $H$-modules $U$ and $V$,  def\/ine a
linear map $C_{U,V}^R$ from $U\otimes V$ to $V\otimes U$ by sending
$(x\otimes y)$ to $ \sum R''y \otimes R'x$ for any $x\in U$, $y \in
V$.

If $V = \oplus_{i =0}^\infty V_i$ is a graded vector space, let
$V_{>n}$ and  $V_{\leq n}$ denote  $ \oplus _{i=n+1} ^\infty  V_i$
and  $ \oplus _{i=0} ^n V_i,$ respectively. We usually denote  $
\oplus _{i=0} ^n V_i$ by $V_{(n)}$. If $\dim V_i <  \infty$
for any natural number~$i$, then~$V$ is called a local f\/inite graded
vector space. We denote by $\iota _i$ the natural injection from~$V_i$ to~$V$ and by $\pi _i$ the corresponding projection from~$V$
to~$V_i$.

Let $H$ be a   bialgebra and a graded coalgebra with an invertible
element $R_n $  in $ H_{(n)}\otimes H_{(n)}$ for any natural $n.$
Assume $R_{n+1} = R_n + W_n$ with $W_n \in H_{(n+1)} \otimes
H_{n+1}+ H_{n+1}\otimes H_{(n+1)}$. $(H, \{R_n \})$ is called a
local quasitriangular bialgebra
 if  $R_n$  is a copairing on $H_{(n)} \otimes H_{(n)}$, and $(H_{(n)}, R_n)$ is almost
cocommutative for any natural number~$n$. In this case, $\{ R_n \}$
is called a~local quasitriangular structure of $H$. Obviously, if
$(H, R)$ is a quasitriangular bialgebra, then $(H,\{ R_n\})$ is a~local quasitriangular bialgebra with $R_0 =R$, $R _i =0$, $H_0 =H$, $H_i
=0 $ for $i >0.$

The following facts are obvious: $\tau ^{-1} =  \tau (id _H \otimes
S )$  (or $ = \tau (S^{-1} \otimes id_A )$) if $A$ is a Hopf
algebra (or  $H$  is a Hopf algebra with invertible antipode) and
$\tau $ is a  skew pairing; $P^{-1} =  (S \otimes id _A)P $  (or  $
= (id_H \otimes S^{-1}) P $) if $H$ is a Hopf algebra (or $A$ is a
Hopf algebra with invertible antipode)  and $P$  is a copairing.

Let $A$ be an algebra and $M$ be an $A$-bimodule. Then the tensor
algebra $T_A(M)$ of $M$ over $A$ is a graded algebra with
$T_A(M)_0=A$, $T_A(M)_1=M$ and $T_A(M)_n=\otimes^n_AM$ for $n>1$.
That is, $T_A(M)=A\oplus(\bigoplus_{n>0}\otimes^n_AM)$ (see
\cite{Ni78}). Let $D$ be another algebra. If $h$ is an algebra map
from $A$ to $D$ and $f$ is an $A$-bimodule map from $M$ to $D$, then
by the universal property of $T_A(M)$ (see \cite[Proposition 1.4.1]{Ni78}) there is a unique algebra map $T_A(h,f): T_A(M)\rightarrow
D$ such that $T_A(h,f)\iota_0=h$ and $T_A(h,f)\iota_1=f$.  One can
easily see that $T_A ( h, f ) = h + \sum _{n>0} \mu ^{n-1}T_n (f )$,
where $T_n(f)$ is the map from $\otimes _A^n M$ to $\otimes_A^nD$
given by $T_n(f)(x_1\otimes x_2 \otimes \cdots \otimes x_n) =
f(x_1)\otimes f(x_2) \otimes \cdots \otimes f(x_n)$, i.e.,
$T_n(f)=f\otimes _A f\otimes_A\cdots\otimes_A f$. Note that $\mu$
can be viewed as a map from $D\otimes _A D$ to $D$.  For the
details, the reader is directed to \cite[Section 1.4]{Ni78}.

Dually, let $C$ be a coalgebra and let $M$ be a $C$-bicomodule. Then
the cotensor coalgebra~$T_C^c(M)$ of~$M$ over $C$ is a graded
coalgebra with $T_C^c(M)_0=C$, $T_C^c(M)_1=M$ and
$T_C^c(M)_n=\Box^n_CM$ for $n>1$. That is,
$T_C^c(M)=C\oplus(\bigoplus_{n>0}\Box^n_CM)$ (see \cite{Ni78}). Let
$D$ be another coalgebra. If $h$ is a coalgebra map from $D$ to $C$
and $f$ is a $C$-bicomodule map from $D$ to $M$ such that $f({\rm
corad}(D))=0$, then by the universal property of $T_C^c(M)$ (see
\cite[Proposition~1.4.2]{Ni78}) there is a unique coalgebra map
$T_C^c(h,f)$ from $D$ to $T_C^c(M)$ such that $\pi_0T_C^c(h,f)=h$
and $\pi_1T_C^c(h,f)=f$. It is not dif\/f\/icult to see that
$T_C^c(h,f)=h+\sum_{n>0}T_n^c(f)\Delta_{n-1}$, where $T_n^c(f)$ is
the map from $\Box_C^n D$ to $\Box_C^n M $ induced by
$T_n(f)(x_1\otimes x_2 \otimes \cdots \otimes x_n)=f(x_1)\otimes
f(x_2) \otimes \cdots \otimes f(x_n)$, i.e., $T^c_n(f)=f\otimes
f\otimes\cdots\otimes f$.

Furthermore, if $B$ is a Hopf algebra and $M$ is a $B$-Hopf
bimodule, then  $T_B(M)$ and $T_B^c (M)$ are two graded Hopf
algebra. Indeed, by \cite[Section 1.4]{Ni78} and \cite[Proposition~1.5.1]{Ni78}, $T_B(M)$ is a graded Hopf algebra with the counit
$\varepsilon=\varepsilon_B\pi_0$ and the comultiplication
$\Delta=(\iota_0\otimes\iota_0)\Delta_B+\sum_{n>0}\mu^{n-1}T_n
(\Delta _M)$, where $\Delta_M=(\iota_0\otimes\iota_1)\delta _M^-+
(\iota_1\otimes\iota_0)\delta_M^+$. Dually,
 $T_B^c(M)$ is a graded Hopf algebra with  multiplication
$\mu = \mu_B (\pi  _0 \otimes \pi  _0) + \sum _{n>0}
 T_n ^c(\mu _M) \Delta _{n-1}$, where $\mu _M = \alpha   _M^-
 (\pi _0 \otimes \pi _1) + \alpha   _M^+(\pi _1 \otimes \pi _0)$.

\section[Yang-Baxter equations]{Yang--Baxter equations} \label{s3}

  Assume  that  $H$ is a   bialgebra and a graded coalgebra with an invertible
element $R_n $  in $ H_{(n)}\otimes H_{(n)}$ for any natural $n$.
For convenience, let (LQT1), (LQT2) and (LQT3) denote (CP1), (CP2)
and (ACO),
 respectively;
\begin{gather*}
{\rm (LQT4)}: \ \ R_{n+1} = R_n + W_n  \quad \mbox{with} \quad W_n \in H_{(n+1)} \otimes
 H_{n+1}+ H_{n+1}\otimes H_{(n+1)};\\
({\rm LQT4}'):   \ \ R_{n+1} = R_n + W_n  \quad \mbox{with} \quad W_n \in H_{n+1} \otimes            H_{n+1}.
\end{gather*}

Then
$(H, \{R_n\})$ is a local quasitriangular bialgebra if and only if
(LQT1), (LQT2), (LQT3) and (LQT4) hold for any natural number $n$.

 Let $H$ be a graded coalgebra
and a bialgebra. A left $H$-module
 $M$ is called an  $H$-module with f\/inite
cycles  if, for any $x\in M,$ there exists a natural number $n_x$
 such that  $H_i x =0$ when $i >n_x$.
 Let $_H{\cal M}^{\rm cf}$ denote the category of all   left  $H$-modules with f\/inite
cycles.

 \begin{lemma}\label{3.1} Let $H$ be a  graded coalgebra and a
bialgebra.
 If $U$ and  $ V $ are  left  $H$-modules with finite cycles,
so is $U \otimes V$.
\end{lemma}

\begin{proof}  For any $x\in U$, $y \in V,$ there exist two natural
numbers $n_x$ and $n_y$, such that $H_{>n_x} x=0 $ and $H_{>n_y} y=0
$. Set $n_{x\otimes y} = 2n_x + 2 n_y$. It is clear that $H_{>n_{x
\otimes y}} (x\otimes y) =0$. Indeed, for any $h \in H_i$ with $i >
n_{x\otimes y}$,
 we see
  \begin{gather*}
 h (x\otimes y) = \sum h_1 x\otimes h_2y
 = 0 \qquad (\mbox{since} \ H   \ \hbox {is graded coalgebra}).\tag*{\qed}
\end{gather*}
\renewcommand{\qed}{}
\end{proof}

\begin{lemma}\label{3.2} Assume that $(H, \{R_n \})$ is a local quasitriangular Hopf algebra.
Then for any   left $H$-modules $U$ and $V$ with finite cycles,
there exists an invertible linear map $C^{\{R_n\}}_{U, V} : U
\otimes V \rightarrow V\otimes U$ such that
\[
C^{\{R_n\}}_{U, V} (x\otimes y) := C^{R_n} (x\otimes
y) = \sum R_n '' y \otimes R_n 'x
\] with $n > 2n _x+ 2n_y$,  for
$x\in U$, $y\in V.$
\end{lemma}

\begin{proof} We f\/irst def\/ine a map $f$ from $U \times V$ to
$V\otimes U$ by sending $(x, y)$ to $\sum R _n'' y \otimes R_n' x $
with $n > 2n_x+2n_y$ for any $x\in U$, $y \in V$. It is clear that $f$
is well def\/ined. Indeed, if $n
> 2n_x +2n_y$, then $C^{R_{n+1}} (x\otimes y) = C^{R_n}(x\otimes y)$
since $R_{n+1} = R_n + W_n$ with $W_n \in H_{(n+1)} \otimes H_{n+1}+
H_{n+1}\otimes H_{(n+1)}$. $f$ is a $k$-balanced function. Indeed,
for $x, y \in U$, $z, w \in V$, $\alpha \in k$, let $n > 2n_x + 2 n_y
+2n_z$. See
\begin{gather*} f(x+y, z) = \sum R_n ''z \otimes R_n ' (x+y)
=\sum R_n ''z \otimes R_n ' x + \sum R_n ''z \otimes R_n ' y\\
\phantom{f(x+y, z)}{} = f(x, z) + f(y, z).
\end{gather*}
  Similarly, we can show that $f(x, z+w) = f (x, z)+ f(x, w)$, $f(x\alpha , z) = f (x, \alpha z)$.
Consequently, there exists a linear map $C^{\{R_n\}}_{U, V} : U
\otimes V \rightarrow V\otimes U$ such that
\[
C^{\{R_n\}}_{U, V} (x\otimes y) = C^{R_n} (x\otimes y)
\] with $n> 2n _x+ 2n_y$,  for $x\in U$, $y\in V.$

 The inverse $(C^{\{R_n\}}_{U, V})^{-1} $ of $C^{\{R_n\}}_{U, V}$ is
 def\/ined by sending  $(y\otimes x)$ to $\sum (R_n^{-1}) 'x \otimes ( R_n^{-1}) '' y$
  with $n > 2n _x+ 2n_y$ for any $x\in U$, $y \in V$.
\end{proof}

\begin{theorem}\label{3.3}
Assume that  $(H, \{R_n \})$ is a local quasitriangular Hopf
algebra. Then $(_H {\cal M}^{\rm cf},$ $ C^{\{R_n\}})$ is a braided tensor category.
\end{theorem}

\begin{proof} Since $H$ is a bialgebra, we have that  $(_H {\cal M},
\otimes, I, a, r, l)$ is a tensor category by \cite[Proposition XI.3.1]{Ka95}. It follows from Lemma~\ref {3.1} that $(_H {\cal M}
^{\rm cf}, \otimes, I, a, r, l)$ is a tensor subcategory of $(_H {\cal
M}, \otimes, I, a, r, l)$. $C^{\{R_n\}}$ is a braiding of $_H {\cal
M}^{\rm cf}$, which can be shown by the way  similar to the proof of
\cite [Proposition VIII.3.1, Proposition XIII.1.4]{Ka95}.
\end{proof}

An $H$-module $M$ is called a graded $H$-module if $M = \oplus _{i
=0}^\infty M_i$ is a graded vector space and $H_i M_j \subseteq
M_{i+j}$ for any natural number $i$ and $j.$

\begin{lemma} \label{3.3'}
Assume that  $H$ is a local finite graded coalgebra and bialgebra.
If $M = \oplus _{i =0}^\infty M_i$ is a~graded $H$-module,  then the
following conditions are equivalent:

(i) $M$ is an $H$-module with finite cycles;

(ii) $Hx$ is finite dimensional for any $x \in M$;

(iii)  $Hx$ is finite dimensional for any homogeneous element $x$ in $M$.
\end{lemma}

\begin{proof} (i) $\Rightarrow $ (ii). For any $x \in M$, there
exists a natural number $n_x$ such that $H_i x =0$ with $i
>n_x$.  Since $H/ (0:x)_H \cong Hx$, where $(0:x)_H := \{h \in H \mid
h \cdot x =0\}$, we have that $Hx$ is f\/inite dimensional.

(ii) $\Rightarrow $ (iii). It is clear.

(iii) $\Rightarrow $ (i). We f\/irst show that, for any homogeneous
element $x \in M_i$, there exists a natural number $n_x$ such that
$H_j x = 0$ with $j >n_x$. In fact, if the above does  not hold,
then there exists $h _j \in H_{n_j} $  such that $h_j x \not= 0$
with $n_1 < n_2 <  \cdots .$ Considering $h_j x \in M_{n_j +i}$ we
have that $\{h_j x \mid j = 1, 2, \dots  \}$
 is linear independent in $ Hx$, which  contradicts to that  $Hx$
is f\/inite dimensional.

 For any  $x \in M$, then  $x = \sum _{i =1}^l x_i$ and $x_i$ is a homogeneous element
 for $i =1, 2,\dots, l$.
There exists  a natural number $n_{x_i}$ such that $H_j x _i = 0 $
with $j >n_{x_i}$. Set $n_x = \sum _{s=1} ^l n_{x_s}$. Then $H_j x
=0$ with $j > n_x$. Consequently, $M$ is an $H$-module with f\/inite
cycles.
\end{proof}

Note that if $M = \oplus _{i =0}^\infty M_i$ is a graded $H$-module,
then both (ii) $\Rightarrow $ (iii) and (iii) $\Rightarrow $ (i)
hold in Lemma~\ref{3.3'}.

Let $_H {\cal M}^{\rm gf}$ and $_H {\cal M}^{\rm gcf}$ denote the category
of all f\/inite dimensional graded  left $H$-modules and the category
of all graded  left $H$-modules with f\/inite cycles. Obviously, they
are two tensor  subcategories  of $_H {\cal M}^{\rm cf}$. Therefore we
have

\begin{theorem} \label{3.4}
Assume that  $(H, \{R_n \})$ is a local quasitriangular Hopf
algebra. Then $(_H {\cal M}^{\rm gf},$ $ C^{\{R_n\}})$ and $(_H {\cal
M}^{\rm gcf}, C^{\{R_n\}})$ are two  braided tensor categories.
\end{theorem}

Therefore, if $M$ is a f\/inite dimensional graded $H$-module (or
$H$-module with f\/inite cycles) over local quasitriangular Hopf
algebra $(H, \{R_n \})$, then $C^{\{R_n\}}_{M, M}$ is a solution of
Yang--Baxter equation on $M$.

It is easy to prove the following.
\begin{theorem} \label {3.5}
Assume that $(H, \{R_n \})$ is a local quasitriangular Hopf
algebra and $R_{n+1} = R_n + W_n$ with $W_n \in H_{n+1} \otimes
H_{(n+1)}$. If $(M, \alpha^-)$ is an $H$-modules   with finite
cycles then $(M, \alpha ^-, \delta ^- )$ is a~Yetter--Drinfeld
$H$-module, where $\delta ^- (x) = \sum R_n'' \otimes R_n' x$ for
any $x\in M$ and $n \geq n_x$.
\end{theorem}

\section{Relation between tensor algebras and co-tensor coalgebras} \label {s1}

\begin{lemma}[See \cite{CM97,DNR01}]\label{1.1}
Let $A$, $B$ and $ C$ be finite dimensional coalgebras, $(M, \delta
_M^-, \delta _M^+)$ and $(N, \delta  _N^-, \delta _N^+)$ be
respectively a finite dimensional $A$-$B$-bicomodule and a finite
dimensional $B$-$C$-bicomodule. Then

(i) $(M^*, \delta  _M^{-}{}^*, \delta  _M^{+}{}^*)$  is a finite
dimensional $A^*$-$B^*$-bimodule;

(ii) $(M \Box _B N, \delta _{M\Box _B N} ^{-},\delta _{M\Box _B N}
^{+})$
 is an $A$-$C$-bicomodule with structure maps $\delta  _{M\Box _B N} ^{-} = \delta  _M^- \otimes id _N$
 and $\delta  _{M\Box _B N} ^{+} = id _M \otimes \delta  _N^+$;

(iii) $M^* \otimes  _{B^*} N^* \cong (M\Box _B N)^*$   (as
$A^*$-$C^*$-bimodules).
\end{lemma}

\begin{lemma}[See \cite{CM97,DNR01}] \label {1.2}   Let $A$, $B$ and $ C$ be finite
dimensional algebras, $(M, \alpha _M^-, \alpha _M^+)$ and $(N,
\alpha  _N^-, \alpha _N^+)$ be
 respectively a finite dimensional $A$-$B$-bimodule and a finite
 dimensional $B$-$C$-bimodule. Then

(i) $(M^*, \alpha _M^{-}{}^*, \alpha _M^{+}{}^*)$ is a finite
dimensional $A^*$-$B^*$-bicomodule;

(ii) $(M \otimes _B N, \alpha _{M\otimes _B N} ^{-},\alpha
_{M\otimes _B N} ^{+})$
 is an $A$-$C$-bimodule with structure maps
 $\alpha _{M\otimes _B N} ^{-} = \alpha  _M ^-\otimes id _N$ and
$\alpha _{M\otimes _B N} ^{+} = id _M \otimes \alpha  _N^+$;

(iii) $M^* \Box _{B^*} N^* \cong (M\otimes _B N)^*$  (as
$A^*$-$C^*$-bicomodules).
\end{lemma}

\begin{proof} (i) and (ii) are easy.

(iii) Consider
\begin{gather*}
M^* \Box _{B^*} N^*  \cong   (M^* \Box _{B^*} N^*)^{**}
 \cong   (M\otimes _B N)^*  \quad (\mbox{by Lemma~\ref{1.1}}).
\end{gather*}
Of course, we can also prove it in the dual way of the proof of
Lemma~\ref{1.1}, by sending $f\otimes_k g $ to $f \otimes_{B^*} g$
for any $f\in M^*$, $g \in N^*$ with $f \otimes g \in M^* \Box _{B^*}
N^*$.
\end{proof}

\begin{theorem}\label {1.3}
If  $A$ is  a finite dimensional algebra and $M$ is a finite
dimensional $A$-bimodule, then  $T_A(M)$ is isomorphic to
subalgebra
 $\sum _{n=0}^\infty (\Box _{A^*}^nM^*)^*$     of $(T_{A^*}
 ^c(M^*))^0$ under map  $\sigma _{T_A(M)}$ and
  $\sigma _{T_A(M)} =  \sigma _A + \sum
_{n>0} \mu ^{n-1} T_n (\sigma _M)$ with  $\mu ^{n-1} T_n (\sigma _M)
= \sigma _{\otimes _A^n M}.$
\end{theorem}

\begin{proof}  We view  $\oplus _{n=0}^\infty (\Box _{A^*}^nM^*)^*$
as inner direct sum of  vector spaces. It is clear that $\sigma _A$
is algebra homomorphism from $A$ to
 $A^{**} \subseteq (T_{A^*} ^c(M^*))^*$ and
 $\sigma _M$ is a $A$-bimodule homomorphism from
$M$ to $M^{**} \subseteq  (T_{A^*} ^c(M^*))^*$. Thus it follows from
\cite[Proposition~1.4.1]{Ni78} that $\phi = \sigma _A + \sum
_{n>0} \mu ^{n-1} T_n (\sigma _M)$ is an algebra homomorphism from
$T_A(M)$ to $(T_{A^*} ^c(M^*))^*$.

It follows from Lemma~\ref {1.1} (iii) that  $\mu ^{n-1}T_n ( \sigma
_M) = \sigma _ {\otimes _A^nM}$. Indeed, we use induction on $n>0$.
Obviously, the conclusion holds when $n=1$. Let  $n>1$,  $N =
\otimes _A^{n-1}M$, $L = (\Box _{A^*}^{n-1}M^*)^*$ and $\zeta = \mu
^{n-2}T_{n-1}(\sigma _M)$. Obviously, $\mu ^{n-1}T_n ( \sigma _M) =
\mu (\zeta \otimes \sigma _M)$. By inductive assumption, $\zeta =
\sigma _N$ is an $A$-bimodule isomorphism from~$N$ to~$L$.
See
\begin{gather*}
\otimes _A^nM = N\otimes _A M  \stackrel {\nu _1} {\cong }  L
\otimes _{A} M^{**} \quad (\mbox{by inductive assumption}) \\
 \phantom{\otimes _A^nM = N\otimes _A M}{} \stackrel {\nu _2} {\cong }  (\Box _{A^*}^{n-1}M^*)^* \otimes
 _{A^{**}}M^{**} \\
 \phantom{\otimes _A^nM = N\otimes _A M}{} \stackrel {\nu _3} {\cong }  ((\Box _{A^*}^{n-1}M^*)\Box
_{A^{*}}M^{*})^*
  \quad (\mbox{by Lemma~\ref{1.1}~(iii)}) \\
 \phantom{\otimes _A^nM = N\otimes _A M}{}= (\Box _{A^*}^n M^*)^*,
\end{gather*}
where $\nu _1 = \sigma _N \otimes _A \sigma _M$, $\nu _2 (f^{**}
\otimes _A g^{**}) = f^{**} \otimes _{A^{**}} g^{**}$  and $\nu _3
(f^{**} \otimes _{A^{**}} g^{**}) = f^{**} \otimes_k  g^{**}$ for
any $f^{**} \in (\Box _{A^*}^{n-1}M^*)^*$, $g^{**} \in M^{**}$.  Now
we have to show $\nu_3\nu_2\nu _1 = \sigma _{\otimes _A ^n M} = \mu
^{n-1} T_n (\sigma _M)$. For any $f^{*} \in \Box _{A^*}^{n-1}M^*$,
$g^{*} \in M^{*}$, $x\in \otimes _A ^{n-1}M$, $y \in M,$ on the one hand
\begin{gather*} \langle \sigma_{\otimes _A^n M} (x\otimes _A y), f^*
\otimes _kg^*\rangle = \langle f^*, x\rangle \langle g^*, y\rangle.
\end{gather*}
On the other hand,
\begin{gather*} \langle \nu _3\nu _2 \nu_1 (x\otimes _A y), f^*
\otimes_k g^*\rangle = \langle \nu _3\nu _2 ( \sigma _N(x)\otimes _{A} \sigma _M(y)), f^* \otimes _kg^*\rangle \\
\phantom{\langle \nu _3\nu _2 \nu_1 (x\otimes _A y), f^* \otimes_k g^*\rangle}{}
= \langle \sigma _M(x)\otimes _{k} \sigma _N(y), f^* \otimes _kg^*\rangle = \langle f^*, x\rangle \langle g^*, y\rangle.
\end{gather*}
Thus $\nu _3\nu _2 \nu_1 = \sigma _{\otimes _A ^n M}$.
 See
\begin{gather*} \langle \mu (\zeta \otimes _A \sigma _M) (x\otimes _A y), f^*
\otimes_k g^*\rangle
=  \langle \zeta (x)\otimes _{A} \sigma _M(y), \Delta ( f^* \otimes _kg^*)\rangle \\
\phantom{\langle \mu (\zeta \otimes _A \sigma _M) (x\otimes _A y), f^* \otimes_k g^*\rangle}{}
= \langle  \zeta (x), f^*\rangle \langle \sigma _M(y), g^*\rangle
=  \langle f^*, x\rangle  \langle g^*, y\rangle.
\end{gather*}
Thus $ \sigma _{\otimes _A ^n M} = \mu ^{n-1} T_n (\sigma _M)$.

Finally, for any $x\in T_A(M)$ with $x = x^{(1)} +x^{(2)} +\cdots
+x^{(n)}$ and $x^{(i)} \in \otimes ^i _A M$,
\begin{gather*}\phi (x) = \sum
_{i=1}^n \phi (x^{(i)})
 = \sum _{i=1}^n \sigma _ {\otimes ^i _A M}(x^{(i)})
= \sigma _{T_A(M)}(x).\tag*{\qed}
\end{gather*}
\renewcommand{\qed}{}
\end{proof}

\section{(Co-)tensor Hopf algebras} \label{s2}

\begin{lemma}\label{2.1} Assume that $B$ is a finite dimensional Hopf algebra
and $M$ is a finite dimensional $B$-Hopf bimodule. Let $A := T_{B}
(M)^{\rm cop}$,
 $H:= T_{B^*}^c(M^*)$.
Then

(i) $ \phi  := \sigma _B + \sum _{i>0} \mu ^{i-1} T_i (\sigma _M)$
is a  Hopf algebra  isomorphism from $T_B(M)$ to the Hopf subalgebra
$\sum _{i=0}^\infty (\Box _{B^*}^iM^*)^*$
     of
$(T_{B^*} ^c(M^*))^0 $;

(ii) Let    $\phi _n : = \phi \mid _{A_{(n)} }$ for any natural
number $n\geq 0$.  Then there exists $\psi_n : (H_{(n)}) ^*
\rightarrow A_{(n)}$ such that $\phi _n \psi _n = id _{(H_{(n)})^*}
$ and $\psi _n \phi _n =id _{A_{(n)}} $, and $\psi _{n+1} (x) = \psi
_n (x)$ for any $x \in (H_{(n)})^*$. Furthermore, $\phi _n$ and
$\psi _n$ preserve the (co)multiplication operations of $T_B(M)$ and
$(T_{B^*}^c (M^*))^0$, respectively.
\end{lemma}

\begin{proof} (i) We f\/irst show that $(\Box _{B^*}^n M^*)^*\subseteq
(T_{B^*} ^c(M^*))^0 $. For any $ f\in (\Box _{B^*}^n M^*)^*$, $\sum
_{i = n+1} ^{\infty}\!\! \Box _{B^*} ^iM^*\!$
 $\subseteq \ker f$ and
$\sum _{i = n+1} ^{\infty} \Box _{B^*} ^iM^*$ is a f\/inite
codimensional ideal of
 $T_{B^*} ^c(M^*).$
Consequently, $f \in (T_{B^*} ^c(M^*))^0$.

Next we show that $\phi  := \sigma _B + \sum _{n>0} \mu ^{n-1} T_n
(\sigma _M)$ (see the proof of Theorem  \ref {1.3}) is a coalgebra
homomorphism from $T_B(M)$ to $ \sum _{i=0}^\infty (\Box _{B^*} ^c
M^*)^*$. For any $x \in \otimes _B^nM$, $f, g \in T_{B^*} ^c(M^*),$ on
the one hand
\begin{gather*} \langle \phi(x), f*g\rangle  =
\langle f*g, x\rangle   \quad  (\mbox{by Theorem~\ref{1.3}}) \\
\phantom{\langle \phi(x), f*g\rangle}{} = \sum _x \langle f, x_1\rangle  \langle g,x_2\rangle
= \sum _x \langle \phi(x_1), f\rangle  \langle \phi(x_2), g\rangle.
\end{gather*}
On the other hand
\begin{gather*} \langle \phi(x), f*g\rangle
 =  \sum \langle (\phi(x))_1, f\rangle \langle (\phi(x))_2, g\rangle,
\end{gather*} since $\phi(x) \in (T_{B^*} ^c(M^*))^0$. Considering
$T_{B^*}^c (M^*) = \oplus _{n\geq 0} \Box _{B^*}^n M^* \cong \oplus
_{n\geq 0} (\otimes _B^nM)^* $ as vector spaces, we have that
$T_{B^*} (M^*)$ is dense in $(T_B(M))^*.$ Consequently,  $\sum
\phi(x_1) \otimes \phi(x_2) = \sum (\phi(x))_1 \otimes (\phi(x))_2$,
i.e.~$\phi$ is a coalgebra homomorphism.

(ii) It follows from Theorem~\ref{1.3}.
\end{proof}

Recall the double cross product $A _\alpha \bowtie _\beta H$,
def\/ined in \cite[p.~36]{ZC01}) and \cite[Def\/inition IX.2.2]{Ka95}.
Assume that $H$ and $A$ are
 two bialgebras;  $(A, \alpha )$ is a left $H$-module coalgebra and
 $(H, \beta )$ is a~right $A$-module coalgebra.
We def\/ine   the multiplication $m_D$, unit $\eta _D$,
comultiplication  $\Delta _D$ and counit $\epsilon _D$  on $A
\otimes H$  as follows:
\begin{gather*}
\mu_D ((a \otimes h) \otimes (b \otimes g) ) = \sum
a\alpha (h_1, b_1) \otimes  \beta (h_2, b_2) g, \\ \Delta _D (a
\otimes h ) = \sum ( a_1 \otimes h_1 \otimes a_2 \otimes h_2),
\end{gather*}
$\epsilon _D = \epsilon _A \otimes \epsilon _H$, $ \eta _D = \eta
_A \otimes \eta _H$ for any $a, b \in A$, $h, g \in H.$ We denote $(A
\otimes H, \mu_D, \eta _D, \Delta _D, \epsilon _D) $
  by
$   A  _{\alpha}  {\bowtie } _{\beta} H$,
 which is  called the double cross product of $A$ and $H$.

\begin{lemma}[See \cite{DT94}]\label {2.2}
Let $H$ and $A$ be two bialgebras. Assume that  $\tau$  is  an
invertible skew pairing on $H \otimes A$. If we
 define
$\alpha (h, a) = \sum \!\tau (h_1, a_1) a_2 \tau ^{-1}(h_2, a_3)$ and
$\beta (h, a) = \sum \!\tau (h_1, a_1) h_2 \tau ^{-1}(h_3, a_2)$
 then  the double cross product $A _\alpha \bowtie _\beta H$
  of $A$ and $H$  is a  bialgebra. Furthermore, if $A$ and $H$ are two
Hopf algebras, then so is  $A _\alpha \bowtie _\beta H$.
\end{lemma}

\begin{proof} We can check that $(A, \alpha )$  is an $H$-module
coalgebra and $(H, \beta )$ is an $A$-module coalgebra step by step.
We can also  check that (M1)--(M4) in \cite[pp.~36--37]{ZC01} hold
step by step. Consequently,  it follows from \cite [Corollary~1.8, Theorem~1.5]{ZC01}  or \cite[Theorem~IX.2.3]{Ka95} that   $A
_\alpha \bowtie _\beta H$ is a Hopf algebra.
\end{proof}

In this case, $A _\alpha \bowtie _\beta H$ can be written as $A
\bowtie _\tau H$.

\begin{lemma}\label {2.3}
Let $H$ and $A$ be two   Hopf algebra. Assume that there exists a Hopf algebra
 monomorphism  $\phi
: A^{\rm cop} \rightarrow H^0$.
Set $ \tau = d_H (\phi \otimes id_H ) C_{H,A}$. Then $ A \bowtie _\tau H$
is  Hopf algebra.
\end{lemma}

\begin{proof} Using \cite[Proposition~2.4]{Ma95a} or the
def\/inition of the evaluation and coevaluation on tensor product,
we can obtain that
 $\tau $ is a skew pairing on $H \otimes A$. Considering Lemma~\ref{2.2},
 we complete the proof.
\end{proof}

\begin{lemma} \label {2.4}
(i) If $H = \oplus_{n=0} ^\infty H_n$ is a graded bialgebra and
$H_0$ has an invertible antipode, then $H$ has an invertible
antipode.

(ii) Assume that $B$ is a finite dimensional Hopf algebra and $M$ is
a
 $B$-Hopf bimodule. Then both
$T_B(M)$ and $T_B^c(M)$ have invertible antipodes.
\end{lemma}

\begin{proof} (i) It is clear that $H^{\rm op}$ is a graded bialgebra
with $(H^{\rm op})_0 = (H_0)^{({\rm op})}$. Thus $H^{\rm op}$ has an antipode by
\cite[Proposition~1.5.1]{Ni78}. However, the antipode of  $H^{\rm op}$
is the inverse of antipode of~$H$.

(ii) It follows from (i).
\end{proof}

\begin{lemma} \label{2.5}
 Let $A = \oplus _{n=0} ^\infty A_n$ and $H = \oplus _{n=0} ^\infty H_n$ be
two graded Hopf algebras with invertible antipodes.
 Let  $\tau$ be a skew pairing on
$(H \otimes A)$ and $P_n$ be a copairing of $H_{(n)} \otimes
A_{(n)}$ for any natural number $n.$  Set $D = A \bowtie _{\tau} H$
and $[P_n] = 1_A \otimes P_n \otimes 1_H.$
 Then $(D, \{ [P_n] \})$  is  almost cocommutative on $D_{(n)}$ if
 and only if
\begin{gather*}
{\rm (ACO1)}: \    \sum P'y_1 \! \otimes P'' \!\otimes y_2 = \! \sum y_4 P' \!\otimes
P''_2 \!\otimes y_2 \tau (y_1, P''_1) \tau ^{-1}(y_3, P''_3) \ \  {\rm  for \ any} \  y \in H_{(n)};\\
{\rm (ACO2)}:  \   \sum x_2\otimes P'\! \otimes x_1P'' = \! \sum x_2 \otimes P'_2\!
\otimes P''x_4\tau (P'_1, x_1) \tau ^{-1}(P'_3, x_3) \ \ {\rm for \ any} \  x\in
A_{(n)}.\!
\end{gather*}
\end{lemma}

\begin{proof} It is clear  that  $(D, \{ [P_n] \})$  is  almost
cocommutative on $D_{(n)}$ if and only if the following holds:
\begin{gather}
  \sum x_2 \otimes y_4P' \otimes x_1 P''_2 \otimes y_2 \tau
(y_1,
P''_1) \tau ^{-1}(y_3, P''_3) \nonumber\\
\qquad{} =  \sum
x_2 \otimes P'_2 y_1\otimes P''x_4 \otimes y_2\tau (P'_1, x_1) \tau
^{-1}(P'_3, x_3)  \label {e2.4''}
\end{gather}
for any $x\in A_{(n)}$, $y\in H_{(n)}$.

Assume that both (ACO1) and (ACO2) hold. See that
\begin{gather*}
  \hbox{the left hand of (\ref {e2.4''})}
 \overset{ \rm by \ (ACO1)}{=}    \sum x_2 \otimes
P'y_1 \otimes x_1 P'' \otimes y_2
\overset{\rm by \ (ACO2)}{=}  \hbox{the right hand of (\ref {e2.4''})}
\end{gather*} for any $x\in A_{(n)}$, $y\in H_{(n)}$.
That is, (\ref {e2.4''}) holds.

Conversely, assume that (\ref {e2.4''}) holds. Thus we have that
\begin{gather*}
  \sum x_2 \otimes P' \otimes x_1 P''_2 \otimes \epsilon _H(1_H )
 \tau (1_H, P''_1) \tau ^{-1}(1_H, P''_3) \\
\qquad{}= \sum
x_2 \otimes P'_2 \otimes P''x_4 \otimes \epsilon _H(1_H)\tau (P'_1,
x_1) \tau ^{-1}(P'_3, x_3)
\end{gather*}
and
\begin{gather*}
  \sum \epsilon _A (1_A) \otimes y_4P' \otimes  P''_2 \otimes y_2
\tau (y_1, P''_1)
 \tau ^{-1}(y_3, P''_3) \\
\qquad{} =  \sum
\epsilon _A(1_A) \otimes P'_2 y_1\otimes P'' \otimes y_2\tau (P'_1,
1_A) \tau ^{-1}(P'_3, 1_A)
\end{gather*} for any $x\in A_{(n)}$, $y\in H_{(n)}$.
Consequently, (ACO1) and (ACO2) hold.
\end{proof}

\begin{lemma} \label{2.6}
 Let $A = \oplus _{n=0} ^\infty A_n$ and $H = \oplus _{n=0} ^\infty H_n$ be
two graded Hopf algebras with invertible antipodes.
 Let  $\tau$ be a skew pairing on
$(H \otimes A)$ and $P_n$ be a copairing of $(H_{(n)} \otimes
A_{(n)})$ with $P_{n+1} = P_n +W_n$ and $W_n \in H _{(n+1)} \otimes
A_{n+1} + H_{n+1} \otimes A_{(n+1)} $ for any natural number $n$.
Set $D = A \bowtie _{\tau} H$. If $\tau (P_n' , x) P_n'' = x$ and
$\tau (y, P_n'' ) P_n' = y$ for any $x \in A_{(n)}$, $y \in H_{(n)}$,
then
  $(D, \{ [P_n] \})$  is a local quasitriangular Hopf algebra.
\end{lemma}

\begin{proof}  It follows from Lemma~\ref{2.2} that $D= \oplus
_{n=0} ^ \infty D_n$ is a  Hopf algebra. Let $D_n = \sum _{i +j =n}
A_i \otimes H_j$. It is clear that $D = \oplus _{n =0} ^\infty D_n$
is a graded coalgebra. We only need to show that $(D, \{ [P_n] \})$
is almost cocommutative on $D_{(n)}$. Now f\/ix $n.$ For convenience,
we denote $P_n$ by $P$ and $Q$ in the following formulae.
For any $x \in A_i$, $y \in H_j$ with $i +j \leq n,$
\begin{gather*}
 \hbox {the right hand  of (ACO1)}
  \overset{\rm by \ (CP1)}{=} \sum y_4 Q' P' \otimes Q_1'' \otimes y_2 \tau (y_1, P'')
\tau
^{-1} (y_3, Q''_2) \\
  \overset{\rm by \ assumption}{=}    \sum y_4Q' y_1 \otimes
Q''_1 \otimes y _2 \tau ^{-1}(y _3, Q''_2)
\overset{\rm by \ (CP1)}{=}   \sum y_4Q' P'y_1 \otimes
P'' \otimes y _2 \tau ^{-1}(y _3, Q'')\\
 =  \sum y_4Q' P'y_1 \otimes
P'' \otimes y _2 \tau (S^{-1}(y _3), Q'')
  \overset{\rm by \ assumption}{=}   \sum y_4 S^{-1}(y_3)
P'y_1 \otimes
P'' \otimes y _2 \\
= \sum  P'y_1 \otimes P'' \otimes y _2
= \hbox{the left hand of (ACO1)}  .
\end{gather*}
 Similarly, we can show that (ACO2) holds on  $A_{(n)}$.
\end{proof}

\begin{theorem} \label {2.7}
Assume that $B$ is a finite dimensional Hopf algebra and $M$ is a
finite dimenaional $B$-Hopf bimodule. Let $A := T_{B} (M)^{\rm cop}$,
 $H:= T_{B^*}^c(M^*)$ and $D = A \bowtie _{\tau}
H$ with $\tau := d_H (\phi  \otimes id )C_{H,A} $. Then $((T_{B}
(M)) ^{\rm cop} \bowtie _\tau T_{B^*} ^c(M^*), \{ R_n\})$ is a local
quasitriangular Hopf algebra. Here $P_n = (id \otimes \psi_n)
b_{H_{(n)} }$,  $R_n = [P_n ] =  1_{B}\otimes (id \otimes \psi_n)
b_{H_{(n)} } \otimes 1_{B^*}$, $\phi $ and $\psi _n$ are defined in
Lemma~{\rm \ref {2.1}}.
\end{theorem}

\begin{proof} By Lemma~\ref {2.4} (ii),  $A$ and $H$ have
 invertible antipodes. Assume that $e_1 ^{(i)}, e_2
^{(i)}, \dots, e_{n_i} ^{(i)}$ is a basis of $H_i$ and $e_1
^{(i)*}, e_2 ^{(i)*}, \dots, e_{n_i} ^{(i)*}$ is an its dual basis
in  $(H_i)^*$. Then $\{e_j ^{(i)} \;|\; i= 0, 1, 2, \dots, n$; $j= 1,
2, \dots, n_i \}$ is a basis of $H_{(n)}$ and $\{e_j ^{(i)*} \; |\;
i=0,  1, 2, \dots, n$; $j = 1, 2, \dots, n_i \}$ is its dual basis
in $(H_{(n)})^*$. Thus  $b_{H_{(n)}} = \sum _ {i =0}^n \sum _{j =1}
^{n_i} e_j^{(i)} \otimes e_j^{(i)*}.$
   See that
\begin {eqnarray*}P_{n+1} &=& \sum
_{i =0} ^{n}\sum _{j =1} ^{n _i}  e_j^{(i)} \otimes \psi _{n+1}
(e_j^{(i)*}) + \sum _{j =1} ^{n _{n+1}}  e_j^{(n+1)}
\otimes  \psi _{n+1}(e_j^{(n+1)*}) \\
&=& P_n +\sum _{j =1} ^{n _{n+1}}  e_j^{(n+1)} \otimes  \psi
_{n+1}(e_j^{(n+1)*}).
\end {eqnarray*}
Obviously, $\sum _{j =1} ^{n _{n+1}}  e_j^{(n+1)} \otimes  \psi
_{n+1}(e_j^{(n+1)*}) \in H_{n+1}\otimes A_{n+1}$. It is clear that
$P_n$ is a copairing on $H _{(n)}\otimes A_{(n)}$ and $\tau$ is a
skew pairing on $H \otimes A$ with  $\sum \tau (P_n' , x) P_n'' = x$
and $\sum \tau (y, P_n'' ) P_n' = y$ for any $x \in A_{(n)}$, $y \in
H_{(n)}$. We complete the proof by Lemma~\ref {2.6}.
\end{proof}

Note that $R_{n+1} = R_n +W_n$ with $W_n \in D_{n+1} \otimes
D_{n+1}$ in the above theorem.

\section{Quiver Hopf algebras} \label{s5}

A quiver $Q=(Q_0,Q_1,s,t)$ is an oriented graph, where  $Q_0$ and
$Q_1$ are the sets of vertices and arrows, respectively; $s$ and $t$
are two maps from  $Q_1$ to $Q_0$. For any arrow $a \in Q_1$, $s(a)$
and $t(a)$ are called its start vertex and end vertex, respectively,
and $a$ is called an arrow from $s(a)$ to $t(a)$. For any $n\geq 0$,
an $n$-path or a path of length $n$ in the quiver $Q$ is an ordered
sequence of arrows $p=a_na_{n-1}\cdots a_1$ with $t(a_i)=s(a_{i+1})$
for all $1\leq i\leq n-1$. Note that a 0-path is exactly a~vertex
and a 1-path is exactly an arrow. In this case, we def\/ine
$s(p)=s(a_1)$, the start vertex of $p$, and $t(p)=t(a_n)$, the end
vertex of $p$. For a 0-path $x$, we have $s(x)=t(x)=x$. Let $Q_n$ be
the set of $n$-paths, $Q_{(n)}$ be the set of $i$-paths with $i \leq
n$ and $Q_{\infty}$ be the set of  all paths in $Q$. Let $^yQ_n^x$
denote the set of all $n$-paths from $x$ to $y$, $x, y\in Q_0$. That
is, $^yQ_n^x=\{p\in Q_n\mid s(p)=x, t(p)=y\}$. A quiver $Q$ is {\it
finite} if $Q_0$ and $Q_1$ are f\/inite sets.

Let $G$ be a group.  Let ${\mathcal K}(G)$ denote the set of
conjugate classes in $G$. $r=\sum_{C\in {\mathcal K}(G)}r_CC$ is
called a {\it ramification} (or {\it ramification data}) of $G$, if
$r_C$ is the cardinal number of a set for any $C\in{\mathcal K}(G)$.
We always assume that the cardinal number of the set $I_C (r)$  is
$r_C$. Let ${\mathcal K}_r(G):=\{C\in{\mathcal K}(G)\mid
r_C\not=0\}=\{C\in{\mathcal K}(G)\mid I_C(r)\not=\varnothing\}$.

Let $G$ be a group. A quiver $Q$ is called {\it a quiver of $G$} if
$Q_0=G$ (i.e., $Q=(G,Q_1,s,t)$). If, in addition, there exists a
ramif\/ication $r$ of $G$ such that the cardinal number of $^yQ_1^x$
is equal to~$r_C$ for any $x, y\in G$ with $x^{-1}y \in C\in
{\mathcal K}(G)$, then $Q$ is called a {\it Hopf quiver with respect
to the ramification data $r$}. In this case, there is a bijection
from $I_C(r)$ to $^yQ_1^x$. Denote by $(Q, G, r)$ the Hopf quiver of
$G$ with respect to $r$.
 $e $ denotes the unit element of $G$. $\{p _g \mid g\in G\}$
denotes the dual basis of $\{g \mid g\in G\}$ of f\/inite group
algebra $kG$.

Let $Q=(G, Q_1, s, t)$ be a quiver of a group $G$. Then $kQ_1$
becomes  a $kG$-bicomodule  under  the natural comodule structures:
\begin{gather}\label{arcom}
\delta^-(a)=t(a)\otimes a,\qquad \delta^+(a)=a\otimes s(a),\qquad a\in
Q_1,
\end{gather} called an {\it arrow comodule}, written as $kQ_1 ^c$.
In this case, the path coalgebra $kQ^c$ is exactly isomorphic to the
cotensor coalgebra $T^c_{kG}(kQ_1^c)$  over $kG$ in a natural way
(see \cite{CM97} and \cite{CR02}). We will regard
$kQ^c=T^c_{kG}(kQ_1^c)$ in the following. Moreover, $kQ_1$ becomes a
$(kG)^*$-bimodule with the module structures def\/ined by
\begin{gather}\label{armod}
p\cdot a:=\langle p, t(a)\rangle a,\qquad a\cdot
p:=\langle p, s(a)\rangle a,\qquad  p\in(kG)^*, \qquad a\in Q_1,
\end{gather}
written as $kQ_1^a$, called an {\it arrow module}. Therefore, we
have a tensor algebra $T_{(kG)^*}(kQ_1^a)$. Note that the tensor
algebra $T_{(kG)^*}(kQ_1^a)$ of $kQ_1^a$ over $(kG)^*$ is exactly
isomorphic to the path algebra~$kQ^a$. We will regard
$kQ^a=T_{(kG)^*}(kQ_1^a)$ in the following.

\begin{lemma}[See \cite{CR02}, Theorem 3.3, and \cite{CR97}, Theorem 3.1] \label {4.1}
Let $Q$ be a quiver over  group~$G$. Then the following statements
are equivalent:

(i) $Q$ is a Hopf quiver.

 (ii) Arrow comodule $kQ_1^c$ admits a
$kG$-Hopf bimodule structure.

If $Q$ is finite, then the above  statements are also equivalent to
the following:

 (iii) Arrow module $kQ_1^a$ admits a
$(kG)^*$-Hopf bimodule structure.

\end{lemma}

Assume that $Q$ is a Hopf quiver. It follows from Lemma \ref {4.1}
that there exist a left $kG$- module structure $\alpha ^-$  and a
right $kG$- module structure $\alpha ^+$ on arrow comodule $(kQ_1^c,
\delta^-, \delta ^+)$ such that  $(kQ_1^c, \alpha ^-, \alpha ^+,
\delta^-, \delta ^+)$ becomes a $kG$-Hopf bimodule, called a
co-arrow Hopf bimodule. We obtain two graded Hopf algebras $T_{kG}
(kQ_1^c)$ and $T_{kG}^c(kQ_1^c)$, called semi-path Hopf algebra and
co-path Hopf algebra, written as $kQ^{s}$ and $kQ^{c}$,
respectively.

Assume that $Q$ is a f\/inite Hopf quiver. Dually, it follows from
Lemma~\ref{4.1} that there exist a left $(kG)^*$-comodule structure
$\delta  ^-$ and a right $(kG)^*$-comodule structure $\delta ^+$ on
arrow module $(kQ_1^a, \alpha^-, \alpha ^+)$ such that  $(kQ_1^a,
\alpha ^-, \alpha ^+, \delta^-, \delta ^+)$ becomes a $(kG)^*$-Hopf
bimodule, called an arrow Hopf bimodule. We obtain two graded Hopf
algebras $T_{(kG)^*} (kQ_1^a)$ and $T_{(kG)^*}^c(kQ_1^a)$, called
path Hopf algebra and semi-co-path Hopf algebra,  written as $kQ^a$
and $kQ^{sc}$, respectively.

From now on, we assume that $Q$ is a f\/inite Hopf quiver on f\/inite
group $G$. Let $\xi _{kQ_1^a}$ denote the linear map from $kQ_1^a$
to $(kQ_1^c)^*$ by sending $a$ to $a^*$ for  any $a\in Q_1$ and $\xi
_{kQ_1^c}$ denote the linear map from $kQ_1^c$ to $(kQ_1^a)^*$ by
sending $a$ to $a^*$ for  any $a\in Q_1$. It is easy to check the
following.

\begin{lemma} \label {4.2} (i) If $(M, \alpha ^-, \alpha ^+, \delta^-, \delta
^+)$ is a finite dimensional $B$-Hopf bimodule and $B$ is a finite
dimensional Hopf algebra, then $(M^*, \delta^{-*}, \delta ^{+*},
\alpha ^{-*}, \alpha ^{+*})$ is a $B^*$-Hopf bimodule.

(ii) If $(kQ_1^c, \alpha ^-, \alpha ^+, \delta^-, \delta ^+)$ is a
co-arrow Hopf bimodule, then there exist  unique left
$(kG)^*$-comodule operation $ \delta _{kQ_1^a}^-$ and right
$(kG)^*$-comodule operation
 $\delta
_{kQ_1^a}^+$ such that $(kQ_1^a, \alpha _{kQ_1^a} ^-, \alpha
_{kQ_1^a}^+,$ \ \ \ $ \delta_{kQ_1^a}^-, \delta_{kQ_1^a} ^+)$
becomes a $(kG)^*$-Hopf bimodule and $\xi_{kQ_1^a}$ becomes a
$(kG)^*$-Hopf bimodule isomorphism from $(kQ_1^a, \alpha _{kQ_1^a}
^-, \alpha _{kQ_1^a}^+, \delta_{kQ_1^a}^-, \delta_{kQ_1^a} ^+)$ to
$((kQ_1^c)^*, \delta^-{}^*, \delta ^+{}^*, \alpha ^-{}^*, \alpha
^+{}^* )$.

(iii) If $(kQ_1^a, \alpha ^-, \alpha ^+, \delta^-, \delta ^+)$ is an
arrow Hopf bimodule, then there exist  unique left $kG$-module
operation $ \alpha _{kQ_1^c}^-$ and right $kG$-module
 $\alpha
_{kQ_1^c}^+$ such that $(kQ_1^c, \alpha _{kQ_1^c} ^-, \alpha
_{kQ_1^c}^+, \delta_{kQ_1^c}^-, \delta_{kQ_1^c} ^+)$ become a
$kG$-Hopf bimodule and $\xi_{kQ_1^c}$ becomes a $kG$-Hopf bimodule
isomorphism from $(kQ_1^c, \alpha _{kQ_1^c} ^-, \alpha _{kQ_1^c}^+, $
$ \delta_{kQ_1^c}^-, \delta_{kQ_1^c} ^+)$ to $((kQ_1^a)^*,
\delta^-{}^*, \delta ^+{}^*, \alpha ^-{}^*, \alpha ^+{}^* )$.

(iv)
 $\xi_{kQ_1^a}$ is a $(kG)^*$-Hopf
 bimodule isomorphism from $\!(kQ_1^a, \alpha _{kQ_1^a} ^-, \alpha
 _{kQ_1^a}^+, \delta_{kQ_1^a}^-, \delta_{kQ_1^a} ^+)\!$ to  $((kQ_1^c)^*,\!$ $
 \delta_{kQ_1^c}^-{}^*, \delta_{kQ_1^c} ^+{}^*, \alpha_{kQ_1^c} ^-{}^*, \alpha _{kQ_1^c}
 ^+{}^* )$ if and only if
 $\xi_{kQ_1^c}$ becomes a $kG$-Hopf bimodule
 isomorphism from $(kQ_1^c, \alpha _{kQ_1^c} ^-, \alpha _{kQ_1^c}^+,
 \delta_{kQ_1^c}^-, \delta_{kQ_1^c} ^+)$ to $((kQ_1^a)^*,
 \delta_{kQ_1^a}^-{}^*, \delta_{kQ_1^a} ^+{}^*, \alpha _{kQ_1^a}^-{}^*, \alpha _{kQ_1^a}^+{}^* )$.
\end{lemma}

Let $B$ be a Hopf algebra and  $^B_B {\cal M}^B_B$ denote the
category of $B$-Hopf bimodules. Let $G$~${\cal H}$opf denote the
category of graded Hopf algebras. Def\/ine $T_B(\psi ) =: T_B (\iota
_0, \iota _1\psi)$ and $T^c_B (\psi) := T_B^c (\pi_0, \psi \pi_1)$
for any $B$-Hopf bimodule homomorphism $\psi.$

\begin{lemma} \label {4.2'} Let $B$  be a Hopf algebra. Then $T_B
$ and $T_B^c$ are two functors from $^B_B {\cal M}^B_B$ to $G$ ${\cal
H }$opf.
\end{lemma}

\begin{proof}
(i) If  $\psi$ is a  $B$-Hopf bimodule homomorphism
from
 $M$ to $M'$, then $T_B(\iota_0, \iota_1\psi)$  is a graded Hopf algebra homomorphism from $T_B(M)$ to
 $T_B(M')$. Indeed, let
$\Phi:=T_B(\iota_0,\iota_1\psi)$. Then both
 $\Delta_{T_{B'}(M')}\Phi$ and $(\Phi\otimes\Phi)\Delta_{T_B(M)}$
 are graded algebra maps from $T_B(M)$ to $T_{B}(M')\otimes
 T_{B}(M')$. Now we show that
 $\Delta_{T_{B}(M')}\Phi=(\Phi\otimes\Phi)\Delta_{T_B(M)}$.
 Considering that $\Phi$  is an algebra homomorphism, we only have to
 show that $\Delta_{T_{B}(M')}\Phi\iota_0=(\Phi\otimes\Phi)\Delta_{T_B(M)}\iota_0$
 and  $\Delta_{T_{B}(M')}\Phi\iota_1=(\Phi\otimes\Phi)\Delta_{T_B(M)}\iota_1$.
Obviously, the f\/irst equation holds. For the second equation, see
\begin{gather*}
  \Delta_{T_{B'}(M')}\Phi\iota_1 = \Delta_{T_{B'}(M')}\iota_1\psi
   = (\iota_0\otimes\iota_1)\delta^-_{M'}\psi+(\iota_1\otimes\iota_0)\delta^+_{M'}\psi\\
 \phantom{\Delta_{T_{B'}(M')}\Phi\iota_1}{}
 = (\iota_0\otimes\iota_1)(id \otimes \psi)\delta^-_M+(\iota_1\otimes\iota_0)(\psi\otimes id)\delta^+_M\\
\phantom{\Delta_{T_{B'}(M')}\Phi\iota_1}{}
 = (\iota_0 \otimes\iota_1\psi)\delta^-_M+(\iota_1\psi\otimes\iota_0)\delta^+_M\\
\phantom{\Delta_{T_{B'}(M')}\Phi\iota_1}{}
= (\Phi\otimes\Phi)(\iota_0\otimes\iota_1)\delta^-_M+(\Phi\otimes\Phi)(\iota_1\otimes\iota_0)\delta^+_M\\
\phantom{\Delta_{T_{B'}(M')}\Phi\iota_1}{}
 = (\Phi\otimes\Phi)[(\iota_0\otimes\iota_1)\delta^-_M+(\iota_1\otimes\iota_0)\delta^+_M]\\
 \phantom{\Delta_{T_{B'}(M')}\Phi\iota_1}{}
   = (\Phi\otimes\Phi)\Delta_{T_B(M)}\iota_1.
  \end{gather*}
  Consequently,  $T_B(\iota_0,\iota_1\psi)$ is a  graded Hopf algebra
  homomorphisms.

(ii) If  $\psi$ is a  $B$-Hopf bimodule   homomorphism from
 $M$ to $M'$,  then $T_B^c(\pi_0, \psi \pi _1)$
 is  a graded Hopf algebra homomorphism from
  $T_B^c(M)$ to $ T_{B}^c(M')$. Indeed,
let $\Psi:=T_{B}^c(\pi _0, \psi\pi_1)$. Then both
$\Psi\mu_{T^c_B(M)}$ and $\mu_{T^c_{B}(M')}(\Psi\otimes\Psi)$ are
graded coalgebra maps from $T^c_B(M)\otimes T^c_B(M)$ to
$T^c_{B}(M')$. Since $T^c_B(M)\otimes T^c_B(M)$ is a graded
coalgebra, ${\rm corad}(T^c_B(M)\otimes
T^c_B(M))\subseteq(T^c_B(M)\otimes
T^c_B(M))_0=\iota_0(B)\otimes\iota_0(B)$. It follows that
$(\pi_1\Psi\mu_{T^c_B(M)})({\rm corad}(T^c_B(M)\otimes
T^c_B(M)))=0$. Thus by the universal property of $T^c_{B}(M')$, in
order to prove
$\Psi\mu_{T^c_B(M)}=\mu_{T^c_{B}(M')}(\Psi\otimes\Psi)$, we only
need to show
$\pi_n\Psi\mu_{T^c_B(M)}=\pi_n\mu_{T^c_{B}(M')}(\Psi\otimes\Psi)$
for $n=0, 1$. However, this follows from a straightforward
computation dual to part~(i). Furthermore, one can see $\Psi(1)=1$.
Hence $\Psi$ is an algebra map, and so a Hopf algebra map.

(iii) It is straightforward to check $T_B (\psi) T_B (\psi') =
T_B(\psi \psi')$ and $T_B ^c(\psi) T_B ^c(\psi') = T_B^c(\psi
\psi')$ for $B$-Hopf bimodule homomorphisms $\psi :   M' \rightarrow
M''$ and $\psi ' :   M \rightarrow M'$.

(iv) $T_B (id _M) = id _{T_B(M)}$ and $T_B ^c(id _M) = id
_{T_B^c(M)}$.
\end{proof}

   \begin{lemma} \label {4.2''} If $\psi$ is a Hopf algebra
   isomorphism from $B$ to $B'$ and $(M, \alpha^-, \alpha ^+, \delta ^-, \delta ^+)$ is \mbox{a~$B$-Hopf}
   bimodule,  then $(M, \alpha ^- (\psi ^{-1}\otimes id_M ), \alpha ^+( id _M\otimes \psi ^{-1}),
   (\psi \otimes id _M) \delta ^-, (id _M \otimes \psi )\delta
   ^+)$ is  \mbox{a~$B'$-Hopf}
   bimodule. Furthermore  $T_B(\iota _0 \psi, \iota_1)$ and $T_B^c(\psi\pi _0,
   \pi_1)$ are graded Hopf algebra isomorphisms from $T_B(M)$ to
   $T_{B'}(M)$ and from $T_B^c(M)$ to
   $T_{B'}^c(M)$, respectively.

  \end{lemma}

 By Lemma \ref {4.2'}
and Lemma \ref {4.2} (ii) and (iii), $T_{(kG)^*} (\iota _0,\iota_1
\xi _{kQ_1^a})$ and $T_{kG} ^c (\pi _0, \xi _{kQ_1^c} \pi_1)$ are
gra\-ded~Hopf algebra isomorphisms from $T_{(kG)^*} (kQ_1^a)$ to
$T_{(kG)^*} ((kQ_1^c)^*)$ and from $T_{kG} ^c(kQ_1^c)$ to
$T_{kG}^c((kQ_1^a)^*)$, respectively. $T_{(kG)^*} ^c (\pi _0, \xi
_{kQ_1^a} \pi_1)$ and $T_{kG}  (\iota _0, \iota _1 \xi _{kQ_1^c})$
are graded Hopf algebra isomorphisms from $T_{(kG)^*} ^c(kQ_1^a)$ to
$T_{(kG)^*} ^c((kQ_1^c)^*)$ and from $T_{kG} (kQ_1^c)$ to
$T_{kG}((kQ_1^a)^*)$, respectively. Furthermore, $(kQ_1^a, kQ_1^c)$,
$(kQ^a, kQ^c)$ and $(kQ^{s}, kQ^{sc})$ are said to be arrow dual
pairings.

\begin{theorem} \label {4.3} Assume that $(Q, G, r)$ is a finite Hopf quiver on finite group $G$. If
$(kQ^a, kQ^c)$ and $(kQ^{s}, kQ^{sc})$ are said to be arrow dual
pairings, then

(i)  $((kQ^a)^{\rm cop} \bowtie _\tau kQ^c, \{R_n\})$ is a local
quasitriangular Hopf algebra. Here{\samepage
\[
R_n = \sum _{g\in G} p_e
\otimes g  \otimes p_g \otimes e + \sum _{q \in Q_{(n)}, \ q \not\in
G} \ \ p_e \otimes q \otimes q \otimes e
\]
 and  $\tau (a,b) =
\delta_{a,b},$   for any two paths~$a$   and~$b$
 in~$Q$, where $\delta_{a,b}$ is the Kronecker symbol.}

(ii) There exist $\tau$ and $\{R_n\}$ such that $((kQ^{s})^{\rm cop}
\bowtie _\tau kQ^{sc}, \{R_n\} )$ becomes  a local quasitriangular
Hopf algebra.
\end{theorem}

\begin{proof}
 (i) Let  $B:= (kG)^*$ and $M:= kQ_1^a. $ Thus  $T_B(M) = kQ^a$.  Since $(kQ^a_1, kQ^c_1)$
 is  an arrow dual pairing, $\xi_{kQ_1^c}$ is a $kG$-Hopf bimodule
 isomorphism from $kQ_1^c$ to $(kQ_1^a)^*$ by Lemma \ref {4.2}.
 See that \begin{gather*}
 T_{B^*} ^c(M^*)  =  T_{(kG)^{**}}^c((kQ_1^a)^*)\\
 \phantom{T_{B^*} ^c(M^*)}{}   \stackrel {\nu_1} {\cong }  T_{kG}^c((kQ_1^a)^*) \qquad (\hbox{by Lemma}~\ref{4.2''})\\
\phantom{T_{B^*} ^c(M^*)}{}   \stackrel {\nu_2} {\cong}  T_{kG}^c(kQ_1^c) \qquad (\hbox{by Lemma}~\ref{4.2} \ \hbox{and  Lemma}~\ref{4.2'}) \\
\phantom{T_{B^*} ^c(M^*)}{} = kQ^c,
 \end{gather*}
 where $\nu _1 = T _{kG}^c ( \sigma ^{-1}_{kG} \pi_0, \pi_1)$, $\nu _2 =
 T _{(kG)^{**}}^c ( \pi_0, (\xi _{kQ_1^c})^{-1}\pi_1)$.

 Let $H = T_{B^*}^c (M^*)$ and $A = T_B (M) ^{\rm cop}.$ By Theorem~\ref{2.7}, $((kQ^a)^{\rm cop} \bowtie _\tau kQ^c, \{R_n\})$ is a local
quasitriangular Hopf algebra. Here $\tau = d_H(\phi \otimes id
)C_{H,A} ((\nu _2 \nu _1)^{-1}\otimes  id _A)$; $R_n =(id _A \otimes
\nu_2\nu_1 \otimes id _A \otimes \nu _2\nu_1 ) (1_{B}\otimes (id
\otimes \psi_n) b_{H_{(n)} } \otimes 1_{B^*} )$, $\phi $ and $\psi
_n$ are def\/ined in Lemma~\ref {2.1}. We have to show that they are
the same as in this theorem. That is,
\begin{gather} \label {4.31e}
d_H(\phi \otimes id )C_{H,A} ((\nu _2 \nu _1)^{-1}\otimes  id
_A)(a\otimes b) =  \delta _{a, b}
\\
(1_{B}\otimes  b_{H_{(n)} } \otimes 1_{B^*} )= (id _A \otimes
(\nu_2\nu_1 )^{-1} \otimes \phi_n \otimes (\nu
_2\nu_1)^{-1} )\nonumber\\
\phantom{(1_{B}\otimes  b_{H_{(n)} } \otimes 1_{B^*} )=}{}
\times \left(\sum _{g\in G} p_e \otimes g  \otimes p_g \otimes e + \sum _{q \in
Q_{(n)}, q \not\in  G}    p_e \otimes q \otimes q \otimes e\right) .\label {4.32e}
 \end{gather}
If $b = b_n b_{n-1} \cdots b_1$ is a $n$-path in $kQ^c$ with $b_i
\in Q_1$ for $i =1, 2, \dots, n,$ then
\begin{gather} \label{4.3e} (\nu_2\nu_1)^{-1}(b) = b_n ^* \otimes b_{n-1}^* \otimes \cdots \otimes
b_1^*  .
 \end{gather}
 If $b \in G$, then $(\nu_2\nu_1)^{-1}(b) = \sigma _{kG} (b)$.
 Consequently, for any $a, b \in
 Q_\infty$ with $b\in kQ^c$ and  $a\in kQ^a$, we have
  \begin{gather*}  d_H(\phi \otimes id )C_{H,A} ((\nu _2 \nu _1)^{-1}\otimes  id _A)
  (a \otimes
  b)
  =  d_H (\phi (a) \otimes (\nu _2 \nu _1)^{-1}(b) ) \\
 \qquad{} =  \left \{ \begin{array}{ll} \langle \phi (a), \sigma _{kG} (b)\rangle  = \langle \sigma _{kG}(b), a\rangle= \delta _{a,b},  &
 \hbox{when} \ \  b \in Q_0,\vspace{1mm}\\
\langle(\nu_2\nu_1)^{-1} , a \rangle \overset{\rm  by \ (\ref{4.3e})}{=} \delta _{a,b}, &
 \hbox{when} \ \ b \not\in Q_0.
 \end{array}\right.
  \end{gather*} Thus (\ref {4.31e}) holds.
Note that $Q_{n} $ is a basis of not only  $(kQ^a)_{n}$ but also
$(kQ^c)_{n}$ for $n>0.$ By (\ref {4.3e}), $\{\phi  (q) \mid q \in
Q_n\}$ is the dual basis of $\{(\nu _2\nu_1) ^{-1} (q) \mid q \in
Q_n\}$ for any $n>0.$ Consequently, (\ref {4.32e}) holds.

(ii) Let  $B:= kG$ and $M:= kQ_1^c. $ Thus  $T_B(M) = kQ^s$. Since
$(kQ^a_1, kQ^c_1)$
 is an arrow dual pairing, $\xi_{kQ_1^a}$ is a $(kG)^*$-Hopf bimodule
 isomorphism from $kQ_1^a$ to $(kQ_1^c)^*$ by Lemma \ref {4.2}. Thus
 \begin{gather*}
 T_{B^*} ^c(M^*)  =  T_{(kG)^*} ^c((kQ_1^c)^*)  \\
\phantom{T_{B^*} ^c(M^*)}{} \stackrel {\nu _3} {\cong }  T_{(kG)^{*}}^c(kQ_1^a) \qquad  (\hbox{by Lemma~\ref{4.2}}~\hbox{and  Lemma~\ref {4.2'}})\\
\phantom{T_{B^*} ^c(M^*)}{}=  kQ^{sc},
 \end{gather*}
 where $\nu _3 =
 T _{(kG)^*}^c ( \pi_0, (\xi _{kQ_1^a})^{-1}\pi_1)$. By Theorem \ref {2.7}, the double cross product of
$T_B(M)^{\rm cop}$ and $T_{B^*}(M^*)$ is a local quasitriangular Hopf
algebra. Consequently, so is the double cross product of
$(kQ^s)^{\rm cop}$ and $kQ^{sc}$.
\end{proof}

Note (LQT4') holds in the above theorem.

\begin{example} \label {4.4}  Let $G= {\bf Z}_2 = (g) = \{e, g\}$ be the group of order
2 with ${\rm char}\, k\not= 2$, $X$ and~$Y$ be respectively the set of
arrows from $g^0$ to $g^0$ and the set of arrows from $g$ to $g$,
and
 $|X| =  |Y| =3 $.
The quiver $Q$ is a Hopf quiver with respect to ramif\/ication $r =
r_{\{e\}} \{e\}$ with $r_{\{e\}} = 3$. Let $\chi ^{(i)} _{e} \in {\bf Z}_2$
 and $a_{y,x} ^{(i)}$ denote the arrow from $x$ to $y$ for $i =1, 2, 3$.
Def\/ine $\delta^- (a_{x,x}^{(i)} ) = x \otimes a_{x,x}^{(i)} $,
 $\delta^+ (a_{x,x}^{(i)} ) = a_{x,x}^{(i)} \otimes x, $
$g \cdot a_{x,x}^{(i)} = a_{gx,gx} ^{(i)}$, $ a_{e,e}^{(i)} \cdot g
= \chi _e ^{(i)} (g)a_{xg,xg} ^{(i)}$ for $x \in G$, $i = 1, 2, 3.$ By
\cite{CR02}, $kQ_1$  is a~$kG$-Hopf bialgebra.  Therefore, it
follows from Theorem \ref {4.3} that
  $((kQ^a)^{\rm cop} \bowtie _\tau kQ^c, \{R_n\})$
is a~local quasitriangular Hopf algebra and for  every f\/inite
dimensional $(kQ^a)^{\rm cop} \bowtie _\tau kQ^c$-module~$M$,
$C^{\{R_n\}}$ is  a solution of Yang--Baxter equations on $M$.
\end {example}

By the way, we obtain the relation between path algebras and path
coalgebras by Theorem~\ref {1.3}.

\begin{corollary} \label {1.4} Let $Q$ be a finite quiver over finite group $G$.
Then Path algebra $kQ^a$ is algebra isomorphic to subalgebra $\sum
_{n=0}^\infty (\Box _{kG}^n kQ^c_1)^*$ of $(kQ^c)^*$.
\end{corollary}

\begin{proof}  Let  $A = (kG)^*$ and $M = kQ_1^c$. It is clear that
$\xi_{kQ_1^c}$ is a $kG$-bicomodule
 isomorphism from $kQ_1^c$ to $(kQ_1^a)^*$.
 See that \begin{gather*}
 T_{B^*} ^c(M^*)  =  T_{(kG)^{**}}^c((kQ_1^a)^*)
   \stackrel {\nu_1} {\cong }  T_{kG}^c((kQ_1^a)^*)
   \stackrel {\nu_2} {\cong}   T_{kG}^c(kQ_1^c)
 =  kQ^c,
 \end{gather*}
 where $\nu _1 = T _{kG}^c ( \sigma ^{-1}_{kG} \pi_0, \pi_1)$, $\nu _2 =
 T _{kG}^c ( \pi_0, (\xi _{kQ_1^c})^{-1}\pi_1)$. Obviously, $\nu _1$
 and $\nu _2$ are coalgebra isomorphism. Consequently, it
follows from Theorem~\ref {1.3} that $kQ^a$ is algebra isomorphism
to subalgebra $\sum _{n=0}^\infty (\Box _{kG}^n kQ_1^c)^*$ of
$(kQ^c)^*$.
\end{proof}

\subsection*{Acknowledgement}  The authors were f\/inancially
supported by the Australian Research Council. S.Z.\ thanks Department
of Mathematics, University of Queensland for hospitality.

\pdfbookmark[1]{References}{ref}
\LastPageEnding

\end{document}